# An efficient 146-line 3D sensitivity analysis code of stress-based topology optimization written in MATLAB

Hao Deng, Praveen S. Vulimiri, and Albert C. To[*]

Department of Mechanical Engineering and Materials Science, University of Pittsburgh, Pittsburgh, PA 15261

*Corresponding author. Email: albertto@pitt.edu

**Abstract**

This paper presents an efficient and compact MATLAB code for three-dimensional stress-based sensitivity analysis. The 146 lines code includes the finite element analysis and p-norm stress sensitivity analysis based on the adjoint method. The 3D sensitivity analysis for p-norm global stress measure is derived and explained in detail accompanied by corresponding MATLAB code. The correctness of the analytical sensitivity is verified by comparison with finite difference approximation. The nonlinear optimization solver is chosen as the Method of moving asymptotes (MMA). Three typical volume-constrained stress minimization problems are presented to verify the effectiveness of sensitivity analysis code. The MATLAB code presented in this paper can be extended to resolve different stress related 3D topology optimization problems. The complete program for sensitivity analysis is given in the Appendix and is intended for educational purposes only.

## 1. Introduction

As computational resources and advanced manufacturing techniques have improved, structural optimization tools such as topology optimization (TO) now contribute a greater role in the design process. TO, first introduced by Bendsoe and Kikuchi [1], determines the ideal values for design variables representing the material distribution in a geometry for an objective, typically for stiffness. TO designs optimized for stiffness can produce complex, organic, and lightweight designs attractive in industrial applications. However, these designs are prone to introduce stress concentrations which could cause the design to fail in use. It is therefore critical to include stress constraints in TO for industrial applications, ensuring the optimal design satisfies all engineering requirements without major modifications deviating from the original optimized result. As noted by Le et al.[2], there are three challenges when implementing stress constraints which have impeded development: the "singularity" phenomenon [3], stress as a local result, and the non-linearity of stress.

As discussed in many papers [4-16], the singularity phenomenon occurs in design variables with low densities still exhibiting a strain value, resulting in an artificially large stress values and complicating the optimization process. Different approaches have been developed to overcome this challenge. Cheng and Guo [6] use a strain relaxation approach to alleviate the strain in these low-density elements. Bruggi [17], Le et al [2], and Holmberg [18] penalize the stress based on the density variable. All three of these works use the Solid Isotropic Material with Penalization (SIMP) method which represents each density variable as a continuous value between 1 and 0. Xia et al. [12] avoids the singularity problem by using the Bi-

Directional Evolutionary Structural Optimization (BESO) method. The BESO method represents each variable discretely as 1 or 0 to represent solid or void material, removing intermediate values and the singularity the intermediate values present. It is also noted the level set approach, another method for TO, has solved stress constrained TO problems as well [4, 9, 10, 19-31]. The level set method also represents material discretely as found in BESO and avoids the singularity phenomenon. Additional advanced methods for stress related topology optimization can be found in Ref [32-35]. This paper, however, will focus on implementing the SIMP method for stress-based topology optimization. The second challenge lies in the fact stress is calculated locally at each element. As noted by Duysinx and Bendsoe [3], this then requires a constraint for each element and increases the necessary computational resources required to compute the problem. Alternatively, Dusyinx and Sigmund [36] introduced a global stress constraint, grouping all local stresses into one measure. However, the global constraint does not provide great control of the local stress behavior. Instead, regional approaches such as those introduced by Paris et al.[37], Le et al.[2], and Holmberg et al.[18] are preferred. These approaches group elemental constraints based on nearby elements or based on the sorted stress values. It was found the regional approaches are a satisfactory compromise, limiting the number of constraints but still providing good local stress control. The third challenge is presented in the nonlinear behavior of stress. Elemental stress measures are greatly affected by the densities in the local neighborhood of the element. As such, a density variable filter [2] is used to smooth the density variables and sensitivities by taking a weighted average of neighboring nodes. This removes the checkerboard effect present in earlier topology optimization algorithms and produces smooth designs and better convergence to the global minimum. Most recently, Senhora et al. [38] presents a consistent topology optimization formulation for mass minimization with local stress constraints by means of the augmented Lagrangian method, named as aggregation-free approach.

In recent years, several computational programs for topology optimization have been published for educational purposes. These codes are helpful for students and engineers to understand the basic mathematical formulation of topology optimization. The detailed reviews for these programs of TO can be found in Ref [39]. Table 1 summarizes the codes published in recent years, including the three major methods: SIMP [40], level set [41], and BESO[42]. This paper aims to educate the reader to implement stress sensitivity analysis of TO. The code provided in this paper is designed to be easy-to-understand and will discuss the implementation for stress-related problems. As stress related TO is a nonlinear optimization, the problems presented are solved using the Method of Moving Asymptotes (MMA) [43].

The remainder of the paper is organized as follows. Section 2 presents the formulation of the p-norm global stress measure. Section 3 describes the three-dimensional finite element method (FEM) formulation and numerical implementation for eight-node hexahedron elements. Section 4 describes the sensitivity analysis and MATLAB implementation. Section 5 formulates in detail the volume-constrained stress minimization problem. Section 6 demonstrates three typical stress related topology optimization problems to verify the effectiveness of sensitivity analysis code, followed by conclusions in Section 7.

**Table 1** Educational computer programs for topology optimization

| Authors and Reference | Programming language | Method |
|---|---|---|
| (2001) O. Sigmund **[44]** | MATLAB | SIMP |
| (2005) Liu Z, Korvink J G et al **[45]** | FEMLAB | Level set |
| (2010) Challis V J **[46]** | MATLAB | Level set |
| (2010) Suresh, Krishnan **[47]** | MATLAB | SIMP (Pareto) |
| (2010) Huang X, Xie Y M. **[42]** | MATLAB | BESO |
| (2011) Andreassen, Erik, et al **[48]** | MATLAB | SIMP |
| (2012) Talischi, Cameron, et al **[49]** | MATLAB | PolyTop |

| (2014) Zegard T, Paulino G H [50] | MATLAB | Ground Structure |
| --- | --- | --- |
| (2015) Aage, Niels, et al. [51] | PETSc | SIMP |
| (2015) Otomori, Masaki, et al. [52] | MATLAB | Level set |
| (2015) Xia L, Breitkopf P [53] | MATLAB | SIMP |
| (2016) Pereira, Anderson, et al [54] | MATLAB | PolyTop |
| (2018) Wei, Peng, et al [55] | MATLAB | Level set |
| (2018) Loyola, Rubén Ansola, et al [56] | MATLAB | SERA |
| (2018) Laurain, Antoine. [57] | FEniCS | Level set |
| (2018) Sanders, Emily D., et al. [58] | MATLAB | PolyMat |
| (2018) Dapogny, Charles, et al. [59] | FreeFem++ | Shape variation |
| (2019) Chen Q, Zhang X, Zhu B. [60] | MATLAB, APDL | SIMP |
| (2019) Gao, Jie, et al. [61] | MATLAB | SIMP |
| (2019) Liang Y, Cheng G. [62] | MATLAB | Integer programming |
| (2020) Smith H, Norato J A.[63] | MATLAB | Geometry projection |
| (2020) Picelli R, et al. [64] | MATLAB | TOBS |
| (2020) Lin H, Xu A, Misra A, et al. [65] | APDL | DER-BESO |
| (2020) Ferrari F, Sigmund O. [66] | MATLAB | SIMP |

## 2. The stress-based topology optimization problem

### 2.1 Stress and Stiffness Penalization

The design domain is discretized with the eight-node hexahedron elements, and each element is assigned with a density variable. The design variable can be written as $x = (x_1, x_2, \cdots x_{nele})$, where $nele$ is total number of elements. The design variable $x$ for the SIMP method is constrained within [0,1], where 0 corresponds to void material and 1 to solid material. To obtain a black-and-white design, a penalization function is introduced to penalize the intermediate density. For stress-based topology optimization problem, the penalization of stress and stiffness for intermediate design variable values is described in Ref [2, 18]. The element stiffness can be expressed based on SIMP penalization function as follows,

$$\boldsymbol{D}_i = x_i^{pl} \boldsymbol{D}_0 \tag{1}$$

where $\boldsymbol{D}_0$ denotes the stiffness of solid material. $pl$ is a penalization factor, which is set to $pl = 3$ in general [48]. The effective stress vector (artificial measurement) is given as follows,

$$\boldsymbol{\sigma}_i = \boldsymbol{D}_0 \boldsymbol{B}_i \boldsymbol{u}_i \tag{2}$$

where $\boldsymbol{u}_i$ denotes the displacement vector at nodes of the $ith$ element, and $\boldsymbol{B}_i$ is strain matrix of $ith$ element. Note that the stress vector $\boldsymbol{\sigma}_i$ in Voigt notation is given as,

$$\boldsymbol{\sigma}_i = \left(\sigma_{ix}, \sigma_{iy}, \sigma_{iz}, \sigma_{ixy}, \sigma_{iyz}, \sigma_{izx}\right)^T \tag{3}$$

The penalized or relaxed stress measure $\widehat{\boldsymbol{\sigma}}_i(x_i)$ is expressed as,

$$\widehat{\boldsymbol{\sigma}}_i(x_i) = \eta(x_i) \boldsymbol{\sigma}_i \tag{4}$$

Several different penalization schemes have been proposed in recent years [2]. A general stress penalization scheme is given by,

$$\eta(x_i) = (x_i)^q \tag{5}$$

where $q$ is a non-negative stress relaxation parameter. As mentioned by Holmberg et al [18], this non-physical penalization scheme is such that the $\hat{\boldsymbol{\sigma}}_i$ equals $\boldsymbol{\sigma}_i$ for solid material and,

$$\lim_{x_i \to 0} \hat{\boldsymbol{\sigma}}_i(x_i) = \mathbf{0} \tag{6}$$

Following this scheme alleviates the singularity phenomenon discussed in the introduction [67]. The definition of von Mises stress can be written as follows,

$$\sigma_{vm,i} = \left(\sigma_{ix}^2 + \sigma_{iy}^2 + \sigma_{iz}^2 - \sigma_{ix}\sigma_{iy} - \sigma_{iy}\sigma_{iz} - \sigma_{iz}\sigma_{ix} + 3\tau_{ixy}^2 + 3\tau_{iyz}^2 + 3\tau_{izx}^2\right)^{\frac{1}{2}} \tag{7}$$

## 2.2 Global p-norm stress measure

The standard p-norm global stress measure is applied for approximating the maximum stress as follows,

$$\sigma_{PN} = \left(\sum_{i=1}^{nele} \hat{\sigma}_{vm,i}{}^p\right)^{1/p} \tag{8}$$

where $\sigma_{vm,i}$ denotes the von Mises stress at the centroid of the $ith$ element, and $p$ is the p-norm aggregation parameter. It is worth to note that the p-norm value approaches the maximum value of $\sigma_{vm}$ when $p \to \infty$, and that

$$\max \hat{\sigma}_{vm} \leq \left(\sum_{i=1}^{nele} \hat{\sigma}_{vm,i}{}^p\right)^{\frac{1}{p}} \tag{9}$$

In general, a greater value of $p$ can provide a more accurate approximation of the maximum von-Mises stress. However, the problem may become ill-conditioned and cause severe oscillations in the optimization process if the value of $p$ is too great ($p > 30$). Thus, an appropriate p-norm value should be selected such that the convergence history is smooth and the maximum stress approximation is sufficient.

## 3. Finite element analysis and MATLAB implementation

The 3D finite element formulation in this paper is based on the code provided by Liu et al[68], where more details regarding efficient numerical implementation are provided. For 3D problems, the constitutive matrix $\boldsymbol{D}_0$ with unit elastic modulus is given by

$$\boldsymbol{D_0} = \frac{1}{(1+v)(1-2v)} \begin{bmatrix} 1-v & v & v & 0 & 0 & 0 \\ v & 1-v & v & 0 & 0 & 0 \\ v & v & 1-v & 0 & 0 & 0 \\ 0 & 0 & 0 & (1-2v)/2 & 0 & 0 \\ 0 & 0 & 0 & 0 & (1-2v)/2 & 0 \\ 0 & 0 & 0 & 0 & 0 & (1-2v)/2 \end{bmatrix} \tag{10}$$

where $v$ is the Poisson's ratio of the isotropic material. Based on finite element method, the stiffness of a linear elastic element for solid material can be formulated as follows,

$$\boldsymbol{K_0} = \int_{-1}^{+1}\int_{-1}^{+1}\int_{-1}^{+1} \boldsymbol{B}^T \boldsymbol{D_0} \boldsymbol{B} d\xi_1 d\xi_2 d\xi_3 \tag{11}$$

where $d\xi_i$ denotes the natural coordinates. Note that the Jacobian matrix is ignored here due to the unit element length. A detailed description can be found in Liu et al [68]. Note that $\boldsymbol{B}$ is strain-displacement matrix, where detailed mathematical formulation and MATLAB implementation for an eight-node hexahedral element can be found in Gao et al [61]. In this paper, the MATLAB implementation of the element stiffness matrix $\boldsymbol{K_0}$ and the strain-displacement matrix $\boldsymbol{B}$ can be done by lines 74-145. The global stiffness $\boldsymbol{K}$ can be obtained by assembling the element stiffness as follows,

$$\boldsymbol{K} = \sum_{i=1}^{nele} \boldsymbol{K}_i(x_i) = \sum_{i=1}^{nele} E(x_i) \boldsymbol{K_0} \tag{12}$$

where $nele$ is total number of elements, and the interpolation coefficient $E(x_i)$ is defined as,

$$E(x_i) = \left(E_m + x_i^{pl}(E_0 - E_m)\right) \tag{13}$$

where $pl$ is a penalty parameter, and $E_0$ and $E_m$ are chosen as: $E_0 = 1, E_m = 1 \times 10^9$. Finally, the nodal displacement from the FEM can be computed through solving the equilibrium equation

$$\boldsymbol{KU} = \boldsymbol{F} \tag{14}$$

where $\boldsymbol{U}$ is the nodal displacement vector, and $\boldsymbol{F}$ is the vector of external loading, which is independent from the design variable. The node ID system is the same as Liu et al. [68], as shown in Fig. 1.

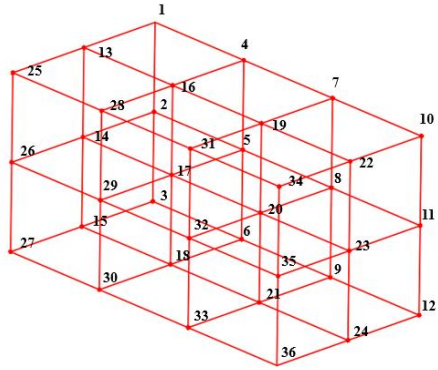

Figure 1. The node ID system of 12 elements

The connectivity matrix edofMat is formulated in lines 20-27, which is used for assembling the global stiffness $\boldsymbol{K}$. The element stiffness matrix and strain-displacement matrix are obtained by line 3. The global stiffness matrix is assembled in lines 28-31 using vectorized MATLAB code. The solution of the nodal displacement vector $\boldsymbol{U}$ is obtained by solving the sparse linear system in line 67. Note that the freedofs and fixeddofs denote the unconstrained and constrained DOFs, respectively. The code example below is for a cantilever beam problem as shown in Fig. 2, where the left side is fixed and an evenly distributed force is applied at the bottom right line. For such boundary conditions, the freedofs and fixeddofs are defined in lines 13,14 and 19. The user can change the boundary and loading conditions by changing the corresponding node IDs.

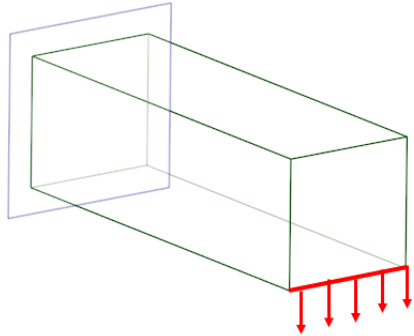

Figure 2. Cantilever beam

## 4. Sensitivity analysis and MATLAB implementation

### 4.1 Sensitivity analysis

To perform the stress-based topology optimization, the sensitivity of the global stress measure with respect to elemental density needs to be provided, a core component of this paper. Note that the total number of elements and nodes are $nele$ and $ndof$, respectively. The derivative of $\sigma_{PN}$ with respect to $jth$ design variable $x_j$ is as follows,

$$\frac{\partial \sigma_{PN}}{\partial x_j} = \sum_{i=1}^{nele} \frac{\partial \sigma_{PN}}{\partial \hat{\sigma}_{vm,i}} \left[ \left(\frac{\partial \hat{\sigma}_{vm,i}}{\partial \hat{\sigma}_i}\right)^T \frac{\partial \eta(x_i)\sigma_i}{\partial x_j} \right] = \sum_{i=1}^{nele} \left( \frac{\partial \sigma_{PN}}{\partial \hat{\sigma}_{vm,i}} \left[ \left(\frac{\partial \hat{\sigma}_{vm,i}}{\partial \hat{\sigma}_i}\right)^T \frac{\partial \eta(x_i)}{\partial x_j} \boldsymbol{\sigma}_i \right] + \frac{\partial \sigma_{PN}}{\partial \hat{\sigma}_{vm,i}} \left[ \left(\frac{\partial \hat{\sigma}_{vm,i}}{\partial \hat{\sigma}_i}\right)^T \eta(x_i) \frac{\partial \boldsymbol{\sigma}_i}{\partial x_j} \right] \right)$$
(15)

The equations above can be rewritten in the following form,

$$\begin{cases} \frac{\partial \sigma_{PN}}{\partial x_j} = T_1 + T_2 \\ T_1 = \sum_{i=1}^{nele} \left( \frac{\partial \sigma_{PN}}{\partial \hat{\sigma}_{vm,i}} \left[ \left(\frac{\partial \hat{\sigma}_{vm,i}}{\partial \hat{\sigma}_i}\right)^T \frac{\partial \eta(x_i)}{\partial x_j} \boldsymbol{\sigma}_i \right] \right) \\ T_2 = \sum_{i=1}^{nele} \left( \frac{\partial \sigma_{PN}}{\partial \hat{\sigma}_{vm,i}} \left[ \left(\frac{\partial \hat{\sigma}_{vm,i}}{\partial \hat{\sigma}_i}\right)^T \eta(x_i) \frac{\partial \boldsymbol{\sigma}_i}{\partial x_j} \right] \right) \end{cases}$$
(16)

The term $\frac{\partial \sigma_{PN}}{\partial \hat{\sigma}_{vm,i}}$ can be expressed as follows,

$$\frac{\partial \sigma_{PN}}{\partial \hat{\sigma}_{vm,i}} = \left( \sum_{i=1}^{nele} (\hat{\sigma}_{vm,i})^p \right)^{\frac{1}{p}-1} (\hat{\sigma}_{vm,i})^{p-1}$$
(17)

Based on the von Mises definition defined in Eq. (7), the derivative of the local element von Mises stress $\hat{\sigma}_{vm,i}$ with respect to the stress vector $\hat{\boldsymbol{\sigma}}_i$ can be written as follows,

$$\begin{cases} \frac{\partial \hat{\sigma}_{vm,i}}{\partial \hat{\sigma}_{ix}} = \frac{1}{2\hat{\sigma}_{vm,i}} (2\hat{\sigma}_{ix} - \hat{\sigma}_{iy} - \hat{\sigma}_{iz}) \\ \frac{\partial \hat{\sigma}_{vm,i}}{\partial \hat{\sigma}_{iy}} = \frac{1}{2\hat{\sigma}_{vm,i}} (2\hat{\sigma}_{iy} - \hat{\sigma}_{ix} - \hat{\sigma}_{iz}) \\ \frac{\partial \hat{\sigma}_{vm,i}}{\partial \hat{\sigma}_{iy}} = \frac{1}{2\hat{\sigma}_{vm,i}} (2\hat{\sigma}_{iz} - \hat{\sigma}_{ix} - \hat{\sigma}_{iy}) \\ \frac{\partial \hat{\sigma}_{vm,i}}{\partial \hat{\tau}_{ixy}} = \frac{3}{\hat{\sigma}_{vm,i}} \hat{\tau}_{ixy} \\ \frac{\partial \hat{\sigma}_{vm,i}}{\partial \hat{\tau}_{ixy}} = \frac{3}{\hat{\sigma}_{vm,i}} \hat{\tau}_{ixy} \\ \frac{\partial \hat{\sigma}_{vm,i}}{\partial \hat{\tau}_{ixy}} = \frac{3}{\hat{\sigma}_{vm,i}} \hat{\tau}_{ixy} \end{cases}$$
(18)

The derivative $\eta(x_i)$ with respect to design variable $x_j$ can be written as,

$$\frac{\partial \eta(x_i)}{\partial x_j} = q x_j^{q-1}$$
(19)

For simplicity of computation, the term $\beta$ is defined as,

$$\beta = (\hat{\sigma}_{vm,i})^{p-1} \left(\frac{\partial \hat{\sigma}_{vm,i}}{\partial \hat{\sigma}_i}\right)^T \frac{\partial \eta(x_i)}{\partial x_j} \boldsymbol{\sigma}_i$$
(20)

The analytical form of term $\frac{\partial \sigma_i}{\partial x_j}$ can be expressed as,

$$\frac{\partial \sigma_i}{\partial x_j} = D_0 B_i \frac{\partial u_i}{\partial x_j} = D_0 B_i L_i \frac{\partial U}{\partial x_j} \qquad (21)$$

where $[L_i]_{24 \times ndof}$ is a 0-1 sparse matrix to extract the nodal displacement of the $ith$ element from the global displacement $[U]_{ndof \times 1}$ ($u_i = L_i U$). Note that material stiffness matrix $[D_0]_{6 \times 6}$ and the elemental strain matrix $[B_i]_{6 \times 24}$ are independent of the design variable $x_j$. Inserting Eq. (21) into the term $T_2$ in Eq. (16), we can find $T_2$ can be rewritten as,

$$T_2 = \sum_{i=1}^{nele} \left( \eta(x_i) \frac{\partial \sigma_{PN}}{\partial \hat{\sigma}_{vm,i}} \left( \frac{\partial \hat{\sigma}_{vm,i}}{\partial \hat{\sigma}_i} \right)^T \frac{\partial \sigma_i}{\partial x_j} \right) = \sum_{i=1}^{nele} \left( \eta(x_i) \frac{\partial \sigma_{PN}}{\partial \hat{\sigma}_{vm,i}} \left( \frac{\partial \hat{\sigma}_{vm,i}}{\partial \hat{\sigma}_i} \right)^T D_0 B_i L_i \frac{\partial U}{\partial x_j} \right) \qquad (22)$$

The adjoint method is applied here to resolve above equation. The term $\frac{\partial U}{\partial x_j}$ can be obtained through differentiating both sides of the equilibrium equation (12) as follows,

$$\frac{\partial K}{\partial x_j} U + K \frac{\partial U}{\partial x_j} = 0 \qquad (23)$$

Therefore, the term $T_2$ can be further written as,

$$T_2 = \sum_{i=1}^{nele} \left( -\eta(x_i) \frac{\partial \sigma_{PN}}{\partial \hat{\sigma}_{vm,i}} \left( \frac{\partial \hat{\sigma}_{vm,i}}{\partial \hat{\sigma}_i} \right)^T D_0 B_i L_i \right) K^{-1} \frac{\partial K}{\partial x_j} U \qquad (24)$$

An adjoint variable $[\lambda]_{ndof \times 1}$ is now defined as,

$$\lambda^T = \left( \sum_{i=1}^{nele} \eta(x_i) \frac{\partial \sigma_{PN}}{\partial \hat{\sigma}_{vm,i}} \left( \frac{\partial \hat{\sigma}_{vm,i}}{\partial \hat{\sigma}_i} \right)^T D_0 B_i L_i \right) K^{-1} \qquad (25)$$

Therefore, adjoint variable $\lambda$ can be calculated from the adjoint equation as follows,

$$K\lambda = \left( \sum_{i=1}^{nele} \eta(x_i) \frac{\partial \sigma_{PN}}{\partial \hat{\sigma}_{vm,i}} (D_0 B_i L_i)^T \frac{\partial \hat{\sigma}_{vm,i}}{\partial \hat{\sigma}_i} \right) \qquad (26)$$

We define the vector $\gamma$ as,

$$\gamma = \left( \sum_{i=1}^{nele} \eta(x_i) \frac{\partial \sigma_{PN}}{\partial \hat{\sigma}_{vm,i}} (D_0 B_i L_i)^T \frac{\partial \hat{\sigma}_{vm,i}}{\partial \hat{\sigma}_i} \right) \qquad (27)$$

Thus, the term $T_2$ can be further simplified to

$$T_2 = -\lambda^T \frac{\partial K}{\partial x_j} U \qquad (28)$$

The sensitivity of the global stiffness matrix $K$ with respect to a design variable $x_j$ equals,

$$\frac{\partial K}{\partial x_j} = \sum_{i=1}^{nele} L_i^T \frac{\partial k_i}{\partial x_j} L_i = L_j^T \left( \frac{\partial k_j}{\partial x_j} \right) L_j = L_j^T \left( pl \cdot x_j^{pl-1} k_j^{(s)} \right) L_j \qquad (29)$$

where $[k_j]_{24 \times 24}$ is the $jth$ elemental stiffness matrix, and $k_j^{(s)}$ is the $jth$ element stiffness with solid material. The detailed expression of element stiffness can be formulated as follows,

$$k_j = x_j^{pl} k_j^{(s)} \text{ and } k_j^{(s)} = \int B_j^T D_0 B_j d\Omega_j \qquad (30)$$

where $\Omega_j$ is the $jth$ element domain. The implementation flow of p-norm stress sensitivity is drawn in Fig. 3.

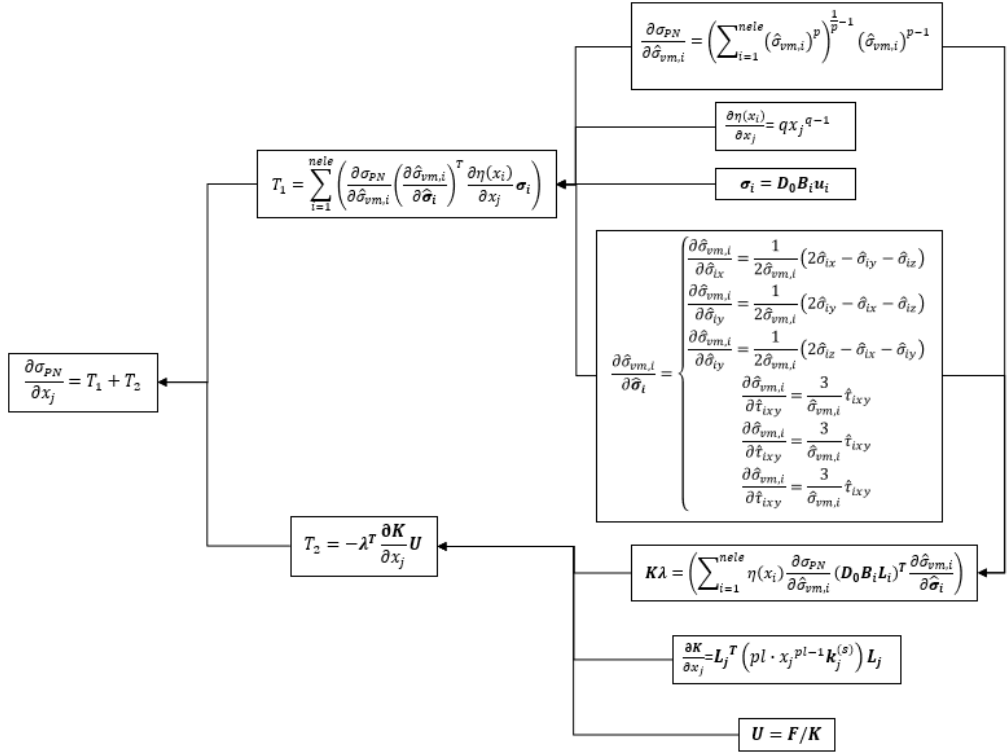

Figure 3. Implementation flow of p-norm stress sensitivity

## 4.2 MATLAB implementation

The MATLAB code for 3D stress sensitivity analysis is executed with the following command:

$$\text{Stress\_3D\_Sensitivity}(x, nelx, nely, nelz, pl, q, p)$$

where nelx, nely and nelz denote the number of elements along $x$, $y$ and $z$ directions. $[x]_{nele \times 1}$ is the elemental density variable constrained in the range of $[0,1]$. $pl$ is the penalization parameter, $q$ is the stress relaxation parameter, and $p$ is the p-norm aggregation parameter. The output of the MATLAB function Stress_3D_Sensitivity is a $nele \times 1$ vector, which is the sensitivity of p-norm stress with respect to the elemental density. The numerical implementation of von Mises stress calculation can be done by lines 33-40 as the following:

```
33 MISES = zeros(nele, 1);
34 S = zeros(nele, 6);
35 for i = 1: nele
36 temp = x(i)^q ∗ (D ∗ B ∗ U(edofMat(i, : )))′;
37 S(i, : ) = temp;
38 MISES(i) = sqrt(0.5 ∗ ((temp(1) − temp(2))^2 + (temp(1) − temp(3))^2. …
39  + (temp(2) − temp(3))^2 + 6 ∗ sum(temp(4: 6). ^2)));
40 end
```

where **S** is a $nele \times 6$ relaxed stress matrix, where each row of matrix **S** denotes a relaxed stress vector $\hat{\boldsymbol{\sigma}}_i$ corresponding to the $ith$ element. MISES is a $nele \times 1$ vector, which stores values of von Mises stress. The objective function (Eq. 8) is formulated in line 46. The numerical implementation of computing $T_1$ can be done as following,

```
47 for i = 1: nele
48 DvmDs(i, 1) = 1/2/MISES(i) * (2 * S(i, 1) − S(i, 2) − S(i, 3));
49 DvmDs(i, 2) = 1/2/MISES(i) * (2 * S(i, 2) − S(i, 1) − S(i, 3));
50 DvmDs(i, 3) = 1/2/MISES(i) * (2 * S(i, 3) − S(i, 1) − S(i, 2));
51 DvmDs(i, 4) = 3/MISES(i) * S(i, 4);
52 DvmDs(i, 5) = 3/MISES(i) * S(i, 5);
53 DvmDs(i, 6) = 3/MISES(i) * S(i, 6);
54 end
55 beta = zeros(nele, 1);
56 for i = 1: nele
57 u = reshape(U(edofMat(i, : ), : )', [], 1);
58 beta(i) = q * (x(i))^(q − 1) * MISES(i)^(p − 1) * DvmDs(i, : ) * D * B * u;
59 end
60 T1 = dpn_dvms * beta;
```

DvmDs is a $nele \times 6$ matrix, where each row represents the derivative of $\hat{\sigma}_{vm,i}$ with respect to $\hat{\boldsymbol{\sigma}}_i$ corresponding to the $ith$ element. beta is a nele × 1 vector, which represents the term $\beta$ in Eq. (20). dpn_dvms is a scaler and given in line 44. It is worth to note that $[U]_{ndof \times 1}$ is global displacement vector and $[u]_{24 \times 1}$ is element displacement vector. The term $T_2$ is computed with the following MATLAB lines:

```
61 gama = zeros(ndof, 1);
62 for i = 1: nele
63 index = index_matrix(: , i);
64 gama(index) = gama(index) + x(i)^q * dpn_dvms * B' * D' * DvmDs(i, : )' * MISES(i).^(p − 1);
65 end
66 lamda = zeros(ndof, 1);
67 lamda(freedofs, : ) = K(freedofs, freedofs)\gama(freedofs, : );
68 T2 = zeros(nele, 1);
69 for i = 1: nele
70 index = index_matrix(: , i);
71 T2(i) = −lamda(index)' * pl * x(i)^(pl − 1) * KE * U(index);
72 end
```

where gama is a ndof × 1 vector defined in Eq. (27). The vector $\boldsymbol{\gamma}$ is computed through lines 62-65 based on the formulation in Eq. (27). The adjoint vector $\boldsymbol{\lambda}$ is computed through lines 66-67, and the term $T_2$ is obtained by lines 68-72. Finally, the sensitivity of p-norm global stress measure with respect to density is obtained in line 73.

## 5. Stress minimization problem formulation

### 5.1 Stress minimization

For simplicity, we focus on the volume-constrained stress minimization problem in this paper, expressed as follows,

$$\begin{cases} \text{minimize } \sigma_{PN} \\ \text{subject to } \int_\Omega d\Omega \leq \bar{V} \end{cases} \tag{31}$$

where $\bar{V}$ is the upper limit for the volume constraint, and $\Omega$ is the design domain. Note that stress sensitivity analysis and corresponding MATLAB code in this paper can be applied to different stress-related problems, such as stress-constrained compliance minimization.

## 5.2 Optimization algorithms

For topology optimization, several advanced optimization algorithms have been proposed in recent years. The Optimality Criteria (OC) [44] is a classical approach for structural optimization problems, followed by recent methods such as Sequential Linear Programming (SLP) and Sequential Quadratic programming (SQP) [68]. Compared to these methods, the Method of moving asymptotes (MMA) proposed by Svanberg [43] is more popular in topology optimization. MMA optimizer approximates the nonlinear optimization problem with the following sub programming problem,

$$\text{find } \tilde{x}$$
$$\text{minimize } -\sum_{i=1}^{n} \left( \frac{x_i^{(k)} - L_i^{(k)}}{\tilde{x}_i - L_i^{(k)}} \right) \frac{\partial \sigma_{PN}}{\partial x_i}(x^{(k)})$$
$$\text{subject to } \tilde{x}^T \mathbf{v} - \bar{v} \leq 0$$
$$\text{where } 0.1 x_i^{(k)} + 0.9 L_i^{(k)} \leq \tilde{x}_i \leq 0.9 U_i^{(k)} + 0.1 x_i^{(k)} \ (i = 1, 2, \cdots n)$$

where $L_i^{(k)}$ and $U_i^{(k)}$ are lower and upper asymptotes and updated in every iteration. The detailed description of updating scheme for $L_i^{(k)}$ and $U_i^{(k)}$ can be found in Svanberg [43]. $x^{(k)}$ is current design point, and $\mathbf{v}$ denotes the volume of element. The MATLAB program (mmasub) implementing the MMA algorithm may be obtained through contacting Prof. Krister Svanberg from KTH. The implementation of the MMA algorithm for the volume-constrained stress minimization problem is straightforward. A sample MATLAB code for computing objective, constraint and their sensitivities for MMA optimizer are as follow,

```
1  function [f0val, df0dx, fval, dfdx] = stress_minimize(x, Hs, H)
2  nelx = 200; nely = 60; nelz = 1; volfrac = 0.3;
3  pl = 3; q = 0.5; p = 10;
4  x(:) = (H ∗ x(:))./Hs;
5  [pnorm, pnorm_sen] = Stress_3D_Sensitivity(x, nelx, nely, nelz, pl, q, p);
6  dv = ones(nely, nelx, nelz)/(nelx ∗ nely ∗ nelz);
7  pnorm_sen(:) = H ∗ (pnorm_sen(:)./Hs);
8  dv(:) = H ∗ (dv(:)./Hs);
9  f0val = pnorm;
10 df0dx = pnorm_sen(:);
11 fval = [mean(x(:)) − volfrac];
12 dfdx = [dv(:)′];
```

Note that Hs and H are filter coefficient vectors; the MATLAB program for computing Hs and H can be found in Liu et al. [68]. volfrac corresponds to the desired volume fraction. f0val is the objective function

(p-norm stress), and vector df0dx is the sensitivity of objective with respect to design variables. Scaler fval is the volume constraint value, and the vector dfdx is the sensitivity of volume constraint with respect to design variables. The filter technique is applied in lines 4, 7 and 8. For density-based topology optimization, filters are critical to avoid numerical instabilities, such as mesh-dependency and checkerboard patterns [69]. Several different filters have been proposed in recent years, and each filter may yield different topology solutions. A review of filter techniques for density-based topology optimization methods can be found in Sigmund [70].

## 6. Numerical example

### 6.1 Sensitivity verification

To verify the correctness of p-norm stress sensitivity, the analytical sensitivity is compared with the sensitivity obtained by the finite difference method. The forward finite difference method for approximating derivatives of a p-norm stress based on a truncated Taylor series expansion can be expressed as,

$$\left[\frac{D\sigma_{PN}(x)}{Dx_i}\right]_f = \frac{\sigma_{PN}(x+\epsilon e_i)-\sigma_{PN}(x)}{\epsilon} \tag{32}$$

where $e_i = [0,0,\cdots,1,\cdots,0,0]^T$ is a unit vector of component $i$ and $\epsilon$ is the perturbation. In general, a smaller $\epsilon$ can get a better approximation for analytical sensitivity, while numerical round-off error will erode the accuracy of the approximation if the value of $\epsilon$ is too small. To quantify the difference between the forward finite difference with the analytical sensitivity, the relative percentage error is introduced as follows,

$$e_f = \left|\frac{\left[\frac{D\sigma_{PN}(x)}{Dx_i}\right]_f - \frac{D\sigma_{PN}(x)}{Dx_i}}{\frac{D\sigma_{PN}(x)}{Dx_i}}\right| \tag{33}$$

where the operator $|\cdot|$ denotes the absolute value. The MATLAB code for sensitivity verification is as follows,

```
1  clear
2  clc
3  nelx = 40;
4  nely = 20;
5  nelz = 1;
6  x = 0.3 * ones(nely, nelx, nelz);
7  x = x(:);
8  d = 1e − 4;
9  e = zeros(nelx * nely, 1);
10 fd = zeros(nelx * nely, 1);
11 pl = 3;
12 q = 0.5;
13 p = 8;
14 for i = 1: nelx * nely
15   x0 = x;
16   x0(i) = x0(i) + d;
17   [pnorm, pnorm_sen] = Stress_3D_Sensitivity(x, nelx, nely, nelz, pl, q, p);
18   [pnorm_f, pnorm_sen] = Stress_3D_Sensitivity(x0, nelx, nely, nelz, pl, q, p);
19   fd(i) = (pnorm_f − pnorm)/d;
20   e(i) = abs((fd(i) − pnorm_sen(i))/pnorm_sen(i));
21 end
```

Due to computational cost, the number of elements in each direction are chosen as nelx = 40, nely = 20, nelz = 1. The perturbation is chosen as $\varepsilon = 1e-4$ in line 8, and the relative percentage error is computed in line 20. The relative percentage error is plotted in Fig. 4 for all design variables $x$. It is worth noting a smaller error can be achieved through adapting the value of perturbation $\varepsilon$. The comparison of sensitivity distribution of analytical and forward finite difference approximation is plotted in Fig. 5, confirming good consistency is achieved. Note that the relative sensitive errors are caused by finite difference step size $\varepsilon$. Adapting step size $\varepsilon$ can obtain a better agreement. We can also verify the sensitivity consistency for any other random density distribution in the range of [0.1,1] by changing the code in line 6 to the following: x = 0.9 * rand(nely, nelx, nelz) + 0.1.

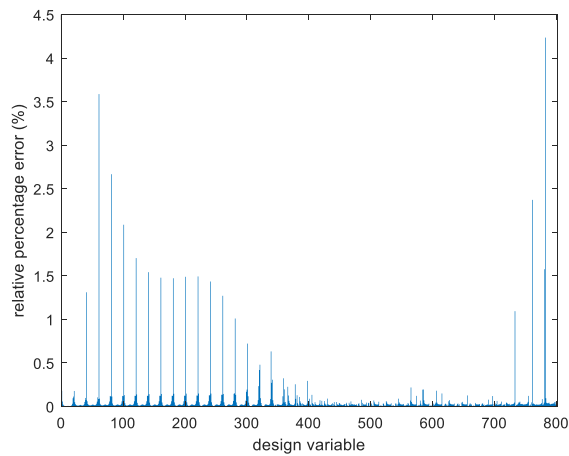

Figure 4. The relative percentage errors

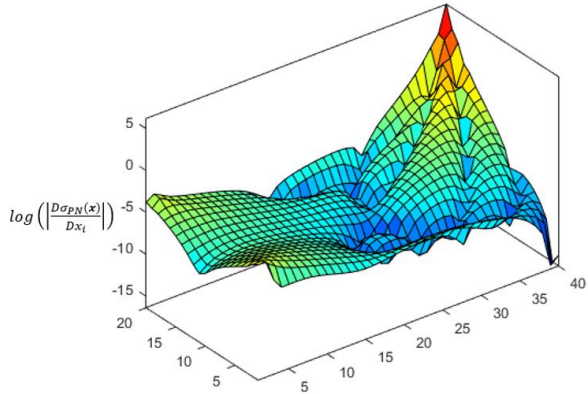 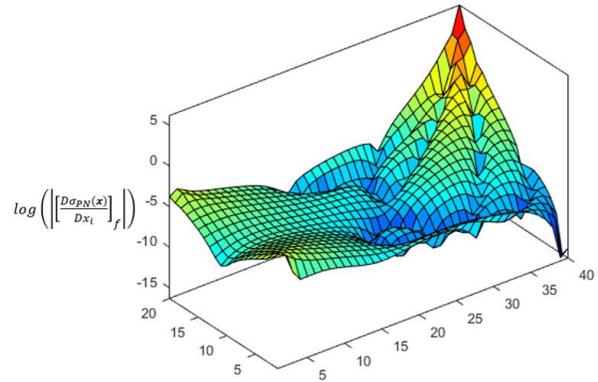

Figure 5. P-norm sensitivity distribution (a) Analytical sensitivity (b) Forward finite difference approximation

### 6.2 Stress minimization for the cantilever beam

The first stress minimization example shown here is a 2D cantilever beam, where $nelx = 200, nely = 60, nelz = 1$ as shown in Fig. 5. The loading is applied on the right-bottom corner, and left side is fully fixed. The volume fraction constraint $\bar{V}$ is chosen as 0.3. The filter radius is $r = 2.5$. The other parameters are selected as: pl=3; q=0.5; p=10. The moving limit of MMA optimizer is set as 0.1. the MMA solver converges after 100 iterations and optimized result is plotted in Fig. 7(a). The von Mises stress distribution is shown in Fig. 7(b). The convergence history is plotted in Fig. 8. The p-norm stress value decreases from initial 33.17 to 3.31 after optimization.

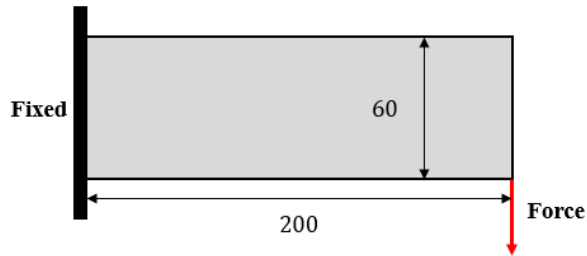

Figure 6. 2D cantilever beam

(a) Optimized material layout            (b) von Mises distribution

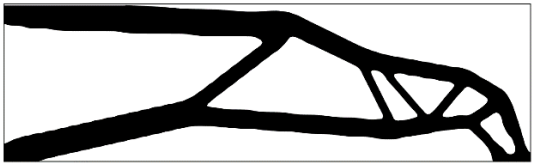 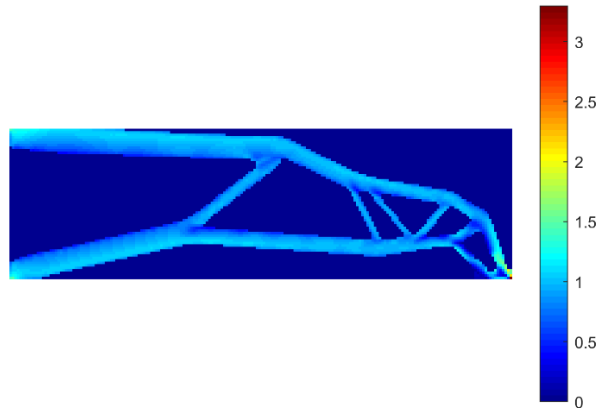

Figure 7. (a) Optimized material layout (b) von Mises distribution

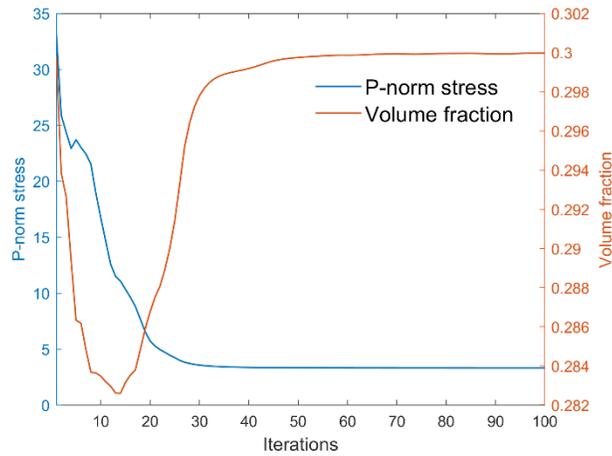

Figure 8. Convergence history

### 6.3 Stress minimization for L-bracket

In this section, a well-known 2D L-bracket benchmark is illustrated in Fig. 9 with the characteristic dimensions. The design domain is discretized by $nelx = 200, nely = 200, nelz = 1$. The passive elements technique is applied as shown in Fig. 9, and the detailed description for passive element implementation can be found in Liu et al [68]. As shown in Fig. 9, the design domain of the L-bracket contains an internal corner which causes a stress singularity. The optimization parameters are chosen as: pl=3; q=0.5; p=10. The filter radius is set to be 2.5. The volume fraction constraint is 0.3. It is worth noting the loading is evenly distributed on three elements at the right corner. Because the loading is uniformed applied on elements, there is no singularity phenomenon nearby loading points. The MATLAB code for boundary and loading can be found and explained in Section 6.4. The optimized result is demonstrated in Fig. 10(a), and von Mises distribution is plotted in Fig. 10(b). The p-norm stress decreases from the initial 56.23 to 5.51 after 120 iterations.

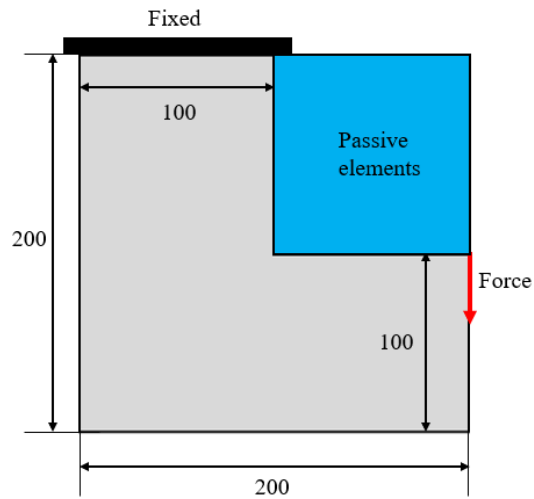

Figure 9. 2D L-bracket example

(a) Optimized material layout          (b) von Mises distribution

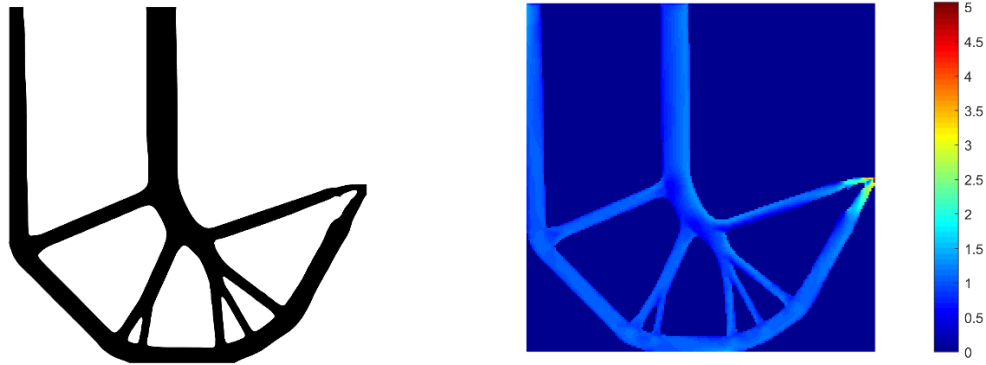

Figure 10. (a) Optimized material layout (b) von Mises distribution

### 6.4 Stress minimization for 3D L-bracket

In this section, we extended the previous 2D L-bracket example to 3D design. The loading and boundary condition are shown in Fig. 11. The design domain is described by $\text{nelx} = 100, \text{nely} = 100, \text{nelz} = 30$. The passive element technique is applied as following MATLAB code,

```
x_t = reshape(x, nely, nelx, nelz);
x_t(1: nelx/2, nely/2: end, : ) = 1e − 4;
x = x_t(: );
```

where x is elemental density vector. The design parameters for 3D L-bracket case is chosen as: pl=3; q=0.5; p=10. The filter radius is 2.5.

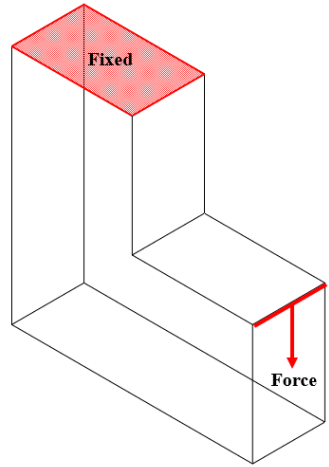

Figure 11. 3D L-bracket example

As shown in Fig. 11, the upper side of L-bracket is fully fixed. The corresponding MATLAB code for such boundary condition is as follows,

```
fixed_node = 1: (nely + 1): ((nely + 1) ∗ (nelx + 1)) − nely;
for j = 1: (nelz + 1)
fixed_node = [fixed_node, 1: (nely + 1): ((nely + 1) ∗ (nelx + 1)) − nely + j ∗ (nely + 1) ∗ (nelx + 1)];
end
fixeddof = [3 ∗ fixed_node − 2,3 ∗ fixed_node − 1,3 ∗ fixed_node]';
```

The vertical loading is applied at the right-upper side, which is uniformly distributed on 3 × nelz elements. The corresponding MATLAB code is shown as following,

```
[il, jl, kl] = meshgrid(nelx, 0, 0: nelz);
loadnid = kl * (nelx + 1) * (nely + 1) + il * (nely + 1) + (nely + 1 − jl);
loaddof1 = (3 * loadnid(:) − 1) − nelx/2 * 3;
loaddof2 = (3 * loadnid(:) − 1) − nelx/2 * 3 + 3;
loaddof3 = (3 * loadnid(:) − 1) − nelx/2 * 3 + 6;
loaddof = union(union(loaddof1, loaddof2), loaddof3);
```

For the 3D L-bracket design problem, the direct solver will take a long time to resolve the large-scale linear systems. The preconditioned conjugate solver built in MATLAB as the function pcg replaces the direct solver. The MATLAB implementation procedure is similar to Ref [68] as follows,

```
tolit = 1e − 8;
maxit = 5000;
M = diag(diag(K(freedofs, freedofs)));
U = pcg(K(freedofs, freedofs), F(freedofs), tolit, maxit, M);
```

where matrix M works as a preconditioner. For computing adjoint vector $\lambda$, we replace the line 67 in MATLAB function Stress_3D_Sensitivity with following code,

```
tolit = 1e − 8;
maxit = 5000;
M = diag(diag(K(freedofs, freedofs)));
lamda(freedofs, : ) = pcg(K(freedofs, freedofs), gama(freedofs), tolit, maxit, M);
```

The optimized material layout for 3D L-bracket example is demonstrated in Fig. 12. Note that the initial sharp corner is replaced by a round shape to avoid a local stress concentration. The p-norm stress decreases significantly from the initial 46.12 to 4.72 after 60 iterations. The von Mises distribution for the optimized result is plotted in Fig. 13. The MATLAB command for plotting von Mises stress is as follows,

$$\text{plot\_von\_Mises(density, von\_Mises)}$$

where the density and von_Mises are both input matrices with dimensions nely × nelx × nelz. The MATLAB code of the function plot_von_Mises can be found in Appendix B. It is worth mentioning the 2D L-bracket example is a special case of the 3D example (nelz = 1).

(a) front view                                            (b) rear view

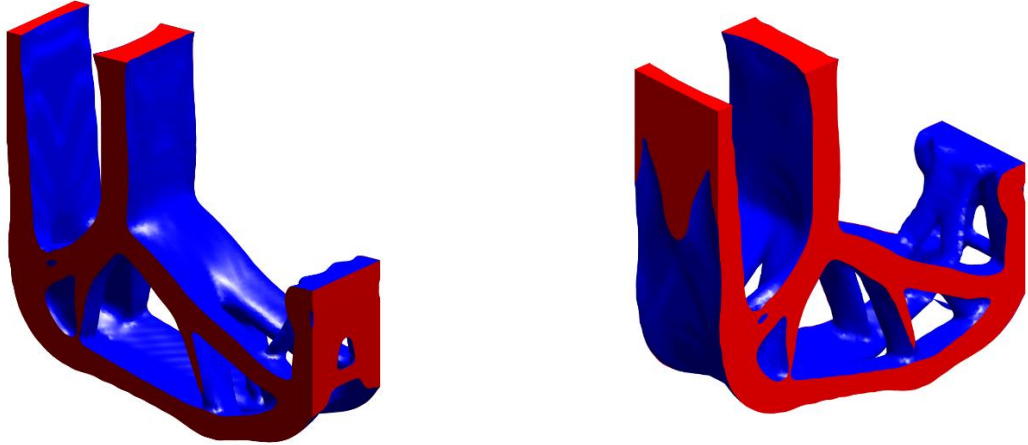

Figure 12. Optimized material layout (a) front view (b) rear view

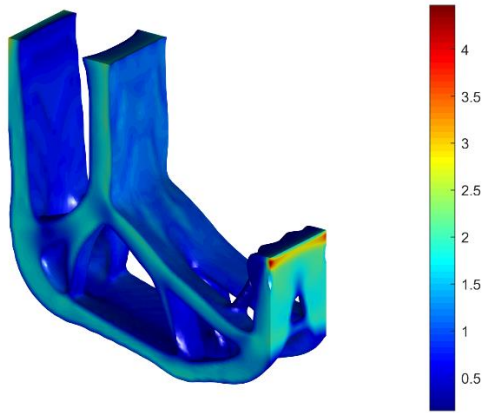

Figure 13. Von Mises stress distribution

## 7. Conclusion

In recent years, stress has become the most important consideration for topology optimization. However, stress-based optimization is challenging due to local nature of computing the stress and highly nonlinear behavior. Local stress is highly dependent on the design, and extreme changes of local stress may occur due to density changes in neighboring regions, meaning stress gradients are highly sensitive to local density changes. Therefore, an accurate, analytical sensitivity is critical for stress related topology optimization problems. In this paper, we derive the sensitivity of the p-norm global stress measure based on a relaxed stress formulation [2], and a corresponding MATLAB code is explained in detail. The finite element formulation in this paper is based on eight-node hexahedral elements, where the design domain is discretized by uniform hexahedral elements. The MMA optimizer is applied for volume-constrained stress minimization problems, where different loading and boundary conditions are explained by MATLAB code. The 146-line code for stress sensitivity analysis is provided for educational purposes, and the users can modify the code based on their requirements. The sensitivity analysis code is written in a compact and vectorized way to enhance the computational efficiency. The complete code is given in the Appendix for users' reference and the source code can be directly download in supplementary material.

## Appendix A: MATLAB Program Stress_3D_Sensitivity

```
1  % A 146-LINE SENSITIVITY ANALYSIS OF STRESS-BASED TOPOLOGY OPTIMIZATION CODE
2  function [pnorm,pnorm_sen]=Stress_3D_Sensitivity(x,nelx,nely,nelz,pl,q,p)
3  [KE,B,D]=brick_stiffnessMatrix();
4  %  MATERIAL PROPERTIES
5  E0 = 1;        % Young's modulus of solid material
6  Emin = 1e-9;   % Young's modulus of void-like material
7  % USER-DEFINED LOAD DOFs
8  [il,jl,kl] = meshgrid(nelx, 0, 0:nelz);          % Coordinates
9  loadnid = kl*(nelx+1)*(nely+1)+il*(nely+1)+(nely+1-jl); % Load Node IDs
10 loaddof = 3*loadnid(:) - 1;              % Load DOFs
11 % USER-DEFINED SUPPORT FIXED DOFs
12 [iif,jf,kf] = meshgrid(0,0:nely,0:nelz);         % Coordinates
13 fixednid = kf*(nelx+1)*(nely+1)+iif*(nely+1)+(nely+1-jf); % Fixed Node IDs
```

```matlab
14  fixeddof = [3*fixednid(:); 3*fixednid(:)-1; 3*fixednid(:)-2]; % Fixed DOFs
15  nele = nelx*nely*nelz;
16  ndof = 3*(nelx+1)*(nely+1)*(nelz+1);
17  F = sparse(loaddof,1,-1,ndof,1);     %External force
18  U = zeros(ndof,1);
19  freedofs = setdiff(1:ndof,fixeddof);
20  nodegrd = reshape(1:(nely+1)*(nelx+1),nely+1,nelx+1);
21  nodeids = reshape(nodegrd(1:end-1,1:end-1),nely*nelx,1);
22  nodeidz = 0:(nely+1)*(nelx+1):(nelz-1)*(nely+1)*(nelx+1);
23  nodeids = repmat(nodeids,size(nodeidz))+repmat(nodeidz,size(nodeids));
24  edofVec = 3*nodeids(:)+1;
25  edofMat = repmat(edofVec,1,24)+ ...
26      repmat([0 1 2 3*nely + [3 4 5 0 1 2] -3 -2 -1 ...
27      3*(nely+1)*(nelx+1)+[0 1 2 3*nely + [3 4 5 0 1 2] -3 -2 -1]],nele,1);
28  iK = reshape(kron(edofMat,ones(24,1))',24*24*nele,1);
29  jK = reshape(kron(edofMat,ones(1,24))',24*24*nele,1);
30  sK = reshape(KE(:)*(Emin+x(:)'.^pl*(E0-Emin)),24*24*nele,1);
31  K = sparse(iK,jK,sK); K = (K+K')/2; %global stiffness matrix assembly
32  U(freedofs,:) = K(freedofs,freedofs)\F(freedofs,:);
33  MISES=zeros(nele,1); %von Mises stress vector
34  S=zeros(nele,6);
35  for i=1:nele
36      temp=x(i)^q*(D*B*U(edofMat(i,:)))';
37      S(i,:)=temp;
38      MISES(i)=sqrt(0.5*((temp(1)-temp(2))^2+(temp(1)-temp(3))^2....
39      +(temp(2)-temp(3))^2+6*sum(temp(4:6).^2)));
40  end
41  nele=size(x(:),1); %total element number
42  [ndof,~]=size(U);
43  DvmDs=zeros(nele,6);
44  dpn_dvms=(sum(MISES.^p))^(1/p-1);
45  index_matrix=edofMat';
46  pnorm=(sum(MISES.^p))^(1/p);
47  for i=1:nele
48      DvmDs(i,1)=1/2/MISES(i)*(2*S(i,1)-S(i,2)-S(i,3));
49      DvmDs(i,2)=1/2/MISES(i)*(2*S(i,2)-S(i,1)-S(i,3));
50      DvmDs(i,3)=1/2/MISES(i)*(2*S(i,3)-S(i,1)-S(i,2));
51      DvmDs(i,4)=3/MISES(i)*S(i,4);
52      DvmDs(i,5)=3/MISES(i)*S(i,5);
53      DvmDs(i,6)=3/MISES(i)*S(i,6);
54  end
55  beta=zeros(nele,1);
56  for i=1:nele
57      u=reshape(U(edofMat(i,:),:)',[],1);
58      beta(i)=q*(x(i))^(q-1)*MISES(i)^(p-1)*DvmDs(i,:)*D*B*u;
59  end
```

```matlab
60  T1=dpn_dvms*beta;
61  gama=zeros(ndof,1);
62  for i=1:nele
63      index=index_matrix(:,i);
64      gama(index)=gama(index)+x(i)^q*dpn_dvms*B'*D'*DvmDs(i,:)'*MISES(i).^(p-1);
65  end
66  lamda=zeros(ndof,1);
67  lamda(freedofs,:)=K(freedofs,freedofs)\gama(freedofs,:);
68  T2=zeros(nele,1);
69  for i=1:nele
70      index=index_matrix(:,i);
71      T2(i)=-lamda(index)'*pl*x(i)^(pl-1)*KE*U(index);
72  end
73  pnorm_sen=T1+T2;
74  function [KE,B,D]=brick_stiffnessMatrix()
75  % elastic matrix formulation
76  nu=0.3;
77  D = 1./((1+nu)*(1-2*nu))*[1-nu nu nu 0 0 0; nu 1-nu nu 0 0 0;...
78      nu nu 1-nu 0 0 0; 0 0 0 (1-2*nu)/2 0 0; 0 0 0 0 (1-2*nu)/2 0;...
79      0 0 0 0 0 (1-2*nu)/2];
80  %stiffness matrix formulation
81  A = [32 6 -8 6 -6 4 3 -6 -10 3 -3 -3 -4 -8;
82      -48 0 0 -24 24 0 0 0 12 -12 0 12 12 12];
83  k = 1/144*A'*[1; nu];
84  K1 = [k(1) k(2) k(2) k(3) k(5) k(5);
85      k(2) k(1) k(2) k(4) k(6) k(7);
86      k(2) k(2) k(1) k(4) k(7) k(6);
87      k(3) k(4) k(4) k(1) k(8) k(8);
88      k(5) k(6) k(7) k(8) k(1) k(2);
89      k(5) k(7) k(6) k(8) k(2) k(1)];
90  K2 = [k(9)  k(8)  k(12) k(6)  k(4)  k(7);
91      k(8)  k(9)  k(12) k(5)  k(3)  k(5);
92      k(10) k(10) k(13) k(7)  k(4)  k(6);
93      k(6)  k(5)  k(11) k(9)  k(2)  k(10);
94      k(4)  k(3)  k(5)  k(2)  k(9)  k(12)
95      k(11) k(4)  k(6)  k(12) k(10) k(13)];
96  K3 = [k(6)  k(7)  k(4)  k(9)  k(12) k(8);
97      k(7)  k(6)  k(4)  k(10) k(13) k(10);
98      k(5)  k(5)  k(3)  k(8)  k(12) k(9);
99      k(9)  k(10) k(2)  k(6)  k(11) k(5);
100     k(12) k(13) k(10) k(11) k(6)  k(4);
101     k(2)  k(12) k(9)  k(4)  k(5)  k(3)];
102 K4 = [k(14) k(11) k(11) k(13) k(10) k(10);
103     k(11) k(14) k(11) k(12) k(9)  k(8);
104     k(11) k(11) k(14) k(12) k(8)  k(9);
105     k(13) k(12) k(12) k(14) k(7)  k(7);
```

```
106      k(10) k(9)  k(8)  k(7)  k(14) k(11);
107      k(10) k(8)  k(9)  k(7)  k(11) k(14)];
108 K5 = [k(1) k(2)  k(8)  k(3)  k(5)  k(4);
109      k(2) k(1)  k(8)  k(4)  k(6)  k(11);
110      k(8) k(8)  k(1)  k(5)  k(11) k(6);
111      k(3) k(4)  k(5)  k(1)  k(8)  k(2);
112      k(5) k(6)  k(11) k(8)  k(1)  k(8);
113      k(4) k(11) k(6)  k(2)  k(8)  k(1)];
114 K6 = [k(14) k(11) k(7)  k(13) k(10) k(12);
115      k(11) k(14) k(7)  k(12) k(9)  k(2);
116      k(7)  k(7)  k(14) k(10) k(2)  k(9);
117      k(13) k(12) k(10) k(14) k(7)  k(11);
118      k(10) k(9)  k(2)  k(7)  k(14) k(7);
119      k(12) k(2)  k(9)  k(11) k(7)  k(14)];
120 KE = 1/((nu+1)*(1-2*nu))*...
121     [ K1  K2  K3  K4;
122      K2'  K5  K6  K3';
123      K3'  K6  K5' K2';
124      K4   K3  K2  K1'];
125 % strain matrix formulation
126 B_1=[-0.044658,0,0,0.044658,0,0,0.16667,0
127 0,-0.044658,0,0,-0.16667,0,0,0.16667
128 0,0,-0.044658,0,0,-0.16667,0,0
129 -0.044658,-0.044658,0,-0.16667,0.044658,0,0.16667,0.16667
130 0,-0.044658,-0.044658,0,-0.16667,-0.16667,0,-0.62201
131 -0.044658,0,-0.044658,-0.16667,0,0.044658,-0.62201,0];
132 B_2=[0,-0.16667,0,0,-0.16667,0,0,0.16667
133 0,0,0.044658,0,0,-0.16667,0,0
134 -0.62201,0,0,-0.16667,0,0,0.044658,0
135 0,0.044658,-0.16667,0,-0.16667,-0.16667,0,-0.62201
136 0.16667,0,-0.16667,0.044658,0,0.044658,-0.16667,0
137 0.16667,-0.16667,0,-0.16667,0.044658,0,-0.16667,0.16667];
138 B_3=[0,0,0.62201,0,0,-0.62201,0,0
139 -0.62201,0,0,0.62201,0,0,0.16667,0
140 0,0.16667,0,0,0.62201,0,0,0.16667
141 0.16667,0,0.62201,0.62201,0,0.16667,-0.62201,0
142 0.16667,-0.62201,0,0.62201,0.62201,0,0.16667,0.16667
143 0,0.16667,0.62201,0,0.62201,0.16667,0,-0.62201];
144 B=[B_1,B_2,B_3];
145 end
146 end
```

## Appendix B: MATLAB Program plot_von_Mises

```matlab
function []=plot_von_Mises(density,von_Mises)
fv=isosurface(density,0.5);
[F1,V1]=isosurface(density,0.5);
[F2,V2]=isocaps(density,0.5);
F3=[F1;F2+length(V1(:,1))];
V3=[V1;V2];
fv.vertices=V3;
fv.faces=F3;
p = patch(fv);
cdata = von_Mises;
isocolors(cdata,p)
p.FaceColor = 'interp';
p.EdgeColor = 'none';
view(3);
axis equal;camlight headlight;lighting phong;
material([0.3,0.6,1,15,1.0]);axis off;colormap('jet');drawnow;
```